\documentclass{amsart}

\usepackage{latexsym}
\usepackage{amscd}
\usepackage{amssymb}
\usepackage{amsmath}
\usepackage{xy}
\xyoption{all}
\usepackage{amsthm}
\usepackage{enumerate}
\usepackage{hyperref}

\newtheorem{thm}{Theorem}[section]
\newtheorem{lem}[thm]{Lemma}
\newtheorem{defn}[thm]{Definition}
\newtheorem{prop}[thm]{Proposition}
\newtheorem{cor}[thm]{Corollary}
\newtheorem{eg}[thm]{Example}
\newtheorem{rmk}[thm]{Remark}

\numberwithin{equation}{section}

\newtheorem*{Thm1}{Theorem 1}
\newtheorem*{Thm2}{Theorem 2}

\newcommand{\leftexp}[2]{{}^{#1}{#2}}

\def\Z{\mathbf{Z}}

\def\P{\mathbf{P}}

\def\ra{\rightarrow}

\def\coh{\operatorname{coh}}
\def\cohc{\operatorname{coh,c}}

\def\Qcoh{\operatorname{QCoh}}
\def\perf{\operatorname{perf}}

\def\Pic{\operatorname{Pic}}

\def\Spec{\operatorname{Spec}}
\def\Func{\operatorname{Func}}

\def\mod{\operatorname{mod}}
\def\Mod{\operatorname{Mod}}

\def\id{\operatorname{id}}

\def\Hom{\operatorname{Hom}}

\def\Aut{\operatorname{Aut}}
\def\vect{\operatorname{vect}}
\def\hocolim{\operatorname{hocolim}}
\def\colim{\operatorname{colim}}
\def\Vect{\operatorname{Vect}}
\def\sing{\operatorname{sing}}

\newcommand{\op}[1]{\operatorname{#1}}

\begin{document}

\author{Matthew Robert Ballard}
\address{Department of Mathematics, University of Pennsylvania,
Philadelphia, PA, USA}
\email{ballardm@math.upenn.edu}

\title[Singular schemes and reconstruction]{Derived categories of sheaves on singular schemes with an application to reconstruction}
\begin{abstract}
 We prove that the bounded derived category of coherent sheaves with proper support is equivalent to the category of locally-finite, cohomological functors on the perfect derived category of a quasi-projective scheme over a field. We introduce the notions of pseudo-adjoints and Rouquier functors and study them. As an application of these ideas and results, we extend the reconstruction result of Bondal and Orlov to Gorenstein projective varieties.
\end{abstract}

\maketitle

\section{Introduction}
 
Understanding schemes through their derived categories of sheaves has become a popular topic in the past couple of decades. However, most results have focused on the case where the scheme is smooth. Let us recall two results most relevant to those in this paper. In \cite{BvB02}, Bondal and van den Bergh proved that the bounded derived category of coherent sheaves, $D^b_{\coh}(X)$, on a smooth and proper variety, $X$, is saturated. Meaning, any covariant or contravariant functor, $\phi: D^b_{\coh}(X) \ra \Vect_k$, satisfying an appropriate boundedness condition, called local-finiteness, is representable. Moreover, $D^b_{\coh}(X)$ is equivalent to either of these categories of functors. For the second result, we take two smooth, projective varieties, $X$ and $Y$, and assume that $\omega_X$ is ample or anti-ample. In \cite{BO01}, Bondal and Orlov proved that, if $X$ and $Y$ have equivalent bounded derived categories, then they must be isomorphic.

In either of these cases, removing the assumption of smoothness sabotages the proofs. For a general projective scheme, $X$, one can define two different categories that reduce to $D^b_{\coh}(X)$ when $X$ is smooth. One is the bounded derived category of coherent sheaves and the other is the smallest triangulated subcategory of $D^b_{\coh}(X)$ containing all finite-rank locally-free sheaves, $D_{\perf}(X)$. The first step in trying to extend the results of Bondal, Orlov, and van den Bergh to a general projective scheme is deciding which category will be the focus of the investigation. The main result of this paper, found in section \ref{sec:lfcohom}, tells us these two categories are very closely related.

\begin{Thm1}
 $D^b_{\coh}(X)$ is equivalent to the category of locally-finite, cohomological functors on $D_{\perf}(X)$.
\label{thm:1}
\end{Thm1}

We prove this result using the machinery of compactly-generated triangulated categories. The functor used in the proof of theorem \ref{thm:1} is a restricted version of the Yoneda functor. Its essential surjectivity was established previously in an appendix to \cite{BvB02}. Fullness follows from the same argument as in \cite{BvB02}. Faithfulness is new.

With our perspective influenced by theorem \ref{thm:1}, we collect some simple ideas and corollaries section \ref{sec:pseudo-ad}. The main application of these ideas is the extension of Bondal and Orlov's result on reconstruction to the case where $X$ and $Y$ are Gorenstein. To realize this extension, we need another idea, a relativization of the notion of a Serre functor. This is presented in section \ref{sec:Roufunc}. In honor of Rouquier's paper \cite{Rou03}, we name this relativization a Rouquier functor. With a Rouquier functor playing the role of the Serre functor and with the ideas and results of the preceding sections, we prove the following in section \ref{sec:recon}.

\begin{Thm2}
 If $X$ is a projective Gorenstein variety with ample or anti-ample canonical bundle and $Y$ is another projective variety with an equivalent perfect derived category (or bounded derived category of coherent sheaves), then $X$ and $Y$ are isomorphic.
\label{thm:2}
\end{Thm2}

This work is a portion of the author's Ph.D. thesis at the University of Washington. The author would like to thank his advisor, Charles F. Doran, for his attention, energy, and suggestions. The author would also like to thank the referee for her/his valuable time and input. While preparing this paper, the author was supported by NSF Research Training Group Grant, DMS 0636606.

\section{Preliminaries}
 
For any category $\mathcal{C}$, the morphism set between objects $A$ and $B$ is denoted as $[A,B]$. All functors are covariant. $k$ is a field.

Before we dive into the bulk of the paper, we will recall some essential ideas and results.

Let $\mathcal{T}$ be a triangulated category possessing all set indexed coproducts. We shall simply say that $\mathcal T$ possesses coproducts. We say an object $C$ of $\mathcal{T}$ is \textbf{compact} (or small) if, for all collections, $X_i$, $i \in I$, of objects in $\mathcal{T}$, the natural map
\begin{displaymath}
 \bigoplus_{i\in I} [C,X_i] \ra [C, \coprod_{i \in I} X_i]
\end{displaymath}
is a isomorphism. The subcategory of compact objects, $\mathcal{T}^c$, is triangulated and closed under taking direct summands.

Given a subcategory $\mathcal{S}$ of $\mathcal{T}$, let $\mathcal{S}^{\perp}$ denote the subcategory of objects, $A$, so that $[S,A]$ is zero for all objects, $S$, of $\mathcal{S}$. We say that $\mathcal{T}$ is \textbf{compactly-generated} if $\left( \mathcal{T}^c \right)^{\perp}$ is zero and the isomorphism classes of objects in $\mathcal T^c$ form a set.

Let $\mathcal R$ be another category possessing coproducts. We say that a functor, $F: \mathcal T \ra \mathcal R$, takes coproducts to coproducts if the natural map $\coprod_{i \in I} F(X_i) \ra F(\coprod_{i \in I} X_i)$ is an isomorphism for all collections of objects, $X_i, i \in I$, from $\mathcal T$. If we assume that $\mathcal R$ has products, we say that $F$ takes coproducts to products if the natural map $F(\coprod_{i \in I} X_i) \ra \prod_{i \in I} F(X_i)$ is an isomorphism for all collections of objects, $X_i, i \in I$.

Let $\mathcal{T}$ be a triangulated category. A functor $H:\mathcal{T}^{\circ} \ra \Mod \ \Z$ is called \textbf{cohomological}, if for each exact triangle
\begin{center}
 \leavevmode
 \begin{xy}
  (-10,10)*+{A}="a"; (10,10)*+{B}="b"; (0,-5)*+{C}="c"; {\ar@{->} "a";"b"}; {\ar@{->} "b";"c"}; {\ar@{-->}^{[1]} "c";"a"}
 \end{xy}
\end{center}
the sequence
\begin{displaymath}
 \cdots \longleftarrow H(C[i-1]) \longleftarrow H(A[i]) \longleftarrow H(B[i]) \longleftarrow H(C[i]) \longleftarrow H(A[i+1]) \longleftarrow \cdots
\end{displaymath}
is exact.
 A functor $H:\mathcal{T} \ra \Mod \ \Z$ is called \textbf{homological}, if for each exact triangle
\begin{center}
 \leavevmode
 \begin{xy}
  (-10,10)*+{A}="a"; (10,10)*+{B}="b"; (0,-5)*+{C}="c"; {\ar@{->} "a";"b"}; {\ar@{->} "b";"c"}; {\ar@{-->}^{[1]} "c";"a"}
 \end{xy}
\end{center}
the sequence
\begin{displaymath}
 \cdots \longrightarrow H(C[i-1]) \longrightarrow H(A[i]) \longrightarrow H(B[i]) \longrightarrow H(C[i]) \longrightarrow H(A[i+1]) \longrightarrow \cdots
\end{displaymath}
is exact.

If $\mathcal{T}$ is compactly-generated, we have a general criteria for representability of functors from $\mathcal{T}$ to abelian groups.

\begin{thm}{(Brown Representability)}
 Let $\mathcal{T}$ be a compactly-generated triangulated category. If $H: \mathcal{T}^{\circ} \ra \Mod \ \Z$ is a functor that takes coproducts to products and is cohomological, then $H$ is representable.
\end{thm}

For a proof, see \cite{Nee96}.

For applications to algebraic geometry, the main example of a compactly-generated triangulated category is the unbounded derived category of quasi-coherent sheaves, $D(X)$, on a quasi-compact and separated scheme, $X$. One can use Brown representability to efficiently study Grothendieck duality for quasi-compact and separated schemes, see \cite{Nee96}.

The category of compact objects of $D(X)$ admits a more geometric characterization. We say that a complex, $E$, in $D(X)$ is \textbf{perfect} if it locally (in the Zariski topology) is quasi-isomorphic to a bounded complex of finite-rank locally-free sheaves. The subcategory of compact objects of $D(X)$ is the subcategory of perfect complexes, see \cite{Nee92}. The subcategory of perfect objects is called the perfect derived category of $X$ and denoted by $D_{\perf}(X)$.

If $X$ possesses the resolution property, i.e. has enough locally-free sheaves, then a complex is perfect if and only if it is quasi-isomorphic to a bounded complex of finite-rank locally-free sheaves, see \cite{TT}. 

We will also need the notion of a \textbf{homotopy colimit} in a triangulated category, $\mathcal T$, possessing coproducts. Given a sequence of morphims in $\mathcal T$,
\begin{displaymath}
 X_0 \overset{f_0}{\ra} X_1 \overset{f_1}{\ra} X_2 \ra \cdots,
\end{displaymath}
The homotopy colimit of $(X_{\bullet},f_{\bullet})$ is denoted by $\hocolim X_i$ and is the unique (up to isomorphism) object fitting into the exact triangle
\begin{center}
 \leavevmode
 \begin{xy}
  (-10,10)*+{\oplus_{i \in \Z}X_i}="a"; (10,10)*+{\oplus_{i \in \Z}X_i}="b"; (0,-5)*+{\hocolim X_i}="c"; {\ar@{->} "a";"b"}; {\ar@{->} "b";"c"}; {\ar@{-->}^{[1]} "c";"a"}
 \end{xy}
\end{center}
With the map $\oplus_{i \in \Z} X_i \ra \oplus_{i \in \Z} X_i$ given by the maps
\begin{center}
 \leavevmode
 \begin{xy}
  (-10,10)*+{X_i}="a"; (10,10)*+{\oplus_{i \in \Z}X_i}="b"; (0,-5)*+{X_i \oplus X_{i+1}}="c"; {\ar@{->} "a";"b"}; {\ar@{->} "c";"b"}; {\ar@{->}_{\op{id}_{X_i} \times -f_i} "a";"c"}
 \end{xy}
\end{center}
for each $i \in \Z$. We will use the following fact in section \ref{sec:pseudo-ad}.
\begin{lem}
 If $C$ is compact object in $\mathcal T$, then $[C,\hocolim X_i]$ is isomorphic to $\colim[C,X_i]$.
\label{lem:hocolim/colim}
\end{lem}
See \cite{Nee92} for a proof. This implies that, given a morphism $\phi: C \ra \hocolim X_i$, there is an $i_0$ and maps $\phi_i: C \ra X_i$ for $i \geq i_0$ so that, for $j > i$, 
\begin{displaymath}
 \phi_j = f_{j-1} \circ \cdots \circ f_i \circ \phi_i 
\end{displaymath}
In addition, for each $i \geq i_0$, morphism $C \overset{\phi_i}{\ra} X_i \ra \hocolim X_i$ is equal to $\phi$. For $D(X)$, we can represent any complex, $E$, as the homotopy colimit of $\tau_{\geq i}E$, the (brutal) truncations of $E$ at the $-i$th step, see \cite{BN}.

\section{Locally-finite cohomological functors on the perfect derived category}
\label{sec:lfcohom}
 
Unless otherwise indicated, $X$ is a quasi-projective scheme over a field $k$ and all categories mentioned are $k$-linear and triangulated. For such an $X$, we denote the bounded derived category of coherent sheaves by $D^b_{\coh}(X)$ and the bounded derived category of coherent sheaves with proper support by $D^b_{\cohc}(X)$. $\Vect_k$ denotes the category of vector spaces over $k$ and $\vect_k$ denotes the category of finite-dimensional vector spaces over $k$.
\begin{defn}
 Let $\mathcal{T}$ be a $k$-linear triangulated category. A $k$-linear functor, $\phi: \mathcal{T} \ra \Vect_k$ (or $\phi: \mathcal{T}^{\circ} \ra \Vect_k$,) is called \textbf{locally-finite} when it satisfies the following condition:
\begin{displaymath}
 \dim_k \left( \bigoplus_{j \in \Z} \phi(A[j]) \right) < \infty \text{ for all } A \in \mathcal{T}.
\end{displaymath}
\end{defn}

Given a $k$-linear triangulated category $\mathcal{T}$, we denote by $\mathcal{T}^{\vee}$ the category of locally-finite cohomological functors on $\mathcal{T}$. We denote by $\leftexp{\vee}{\mathcal{T}}$ the category of locally-finite homological functors on $\mathcal{T}$.

The main result of this section is the following.
\begin{thm}
 Let $X$ be a quasi-projective scheme over a field $k$. Then the restricted Yoneda functor is an equivalence between $D^b_{\cohc}(X)$ and $D_{\perf}(X)^{\vee}$.
\label{thm:lfcohom}
\end{thm}

The proof will be accomplished through the following series of lemmas. The restricted Yoneda functor will be defined in the process of the proof.

\begin{lem}
 Let $\mathcal{T}$ be a compactly-generated triangulated category. Any cohomological functor $F: (\mathcal{T}^c)^{\circ} \ra \vect_k$ is representable by an object of $\mathcal{T}$. Any natural transformation between such functors is induced by a morphism of their representing objects.
\label{lem:rep}
\end{lem}

\proof The statement on objects is lemma 2.14 of \cite{CKN01}. We recall the proof for the sake of completeness. Given a cohomological functor, $F:(\mathcal{T}^c)^{\circ} \ra \Vect_k$, we let $D$ denote the dualization functor on $\Vect_k$, i.e. $D(V) := \op{Hom}_k(V,k)$. Let us set $G = D \circ F$. $G$ is homological. By proposition 2.3 of \cite{Kra}, there is a unique extension of $G$ to a functor, $G^e: \mathcal T \to \Vect_k$, that is homological and preserves coproducts. We now apply $D$ again and use Brown representability to deduce that $D \circ G^e$ is represented by an object $X$ of $\mathcal{T}$. The restriction of $D \circ G$ to $\mathcal{T}^c$ is isomorphic to $D \circ D \circ F$. Since $F$ lands in $\vect_k$, $D^2$ cancels out and $X$ is the object we seek.

Any natural transformation, $\nu: F_1 \ra F_2$, between two cohomological functors induces a natural transformation $D \circ \nu: D \circ F_2 \ra D \circ F_1$. Let us denote $D \circ F_1$ by $G_1$ and $D \circ F_2$ by $G_2$. The construction of lemma 2.2 of \cite{Kra} yields a natural tranformation $(D \circ \nu)^e :G_2^e \to G_1^e$ extending $D \circ \nu$.  Applying $D$ again, we get a natural transformation: $D \circ (D \circ \nu)^e: D \circ G_1 \ra D \circ G_2$. $D \circ G_1$ and $D \circ G_2$ are represented by $X_1$ and $X_2$, respectively. By the Yoneda lemma, $D \circ (D \circ \nu)^e$ is represented by a morphism $\phi: X_1 \ra X_2$. The restriction of $D \circ (D \circ \nu)^e: D \circ G_1 \ra D \circ G_2$ to $\mathcal T^c$ is isomorphic to $D \circ D \circ \nu: D \circ D \circ F_1 \ra D \circ D \circ F_2$ which is isomorphic to $\nu: F_1 \ra F_2$ as each $F_i(C)$ is finite dimensional when $C$ is compact. Thus, $\phi: X_1 \ra X_2$ induces $\nu: F_1 \ra F_2$.\qed


We can rephrase this as follows. The inclusion $\mathcal{T}^c \hookrightarrow \mathcal{T}$ induces a restricted Yoneda functor.
\begin{gather*}
 \mathcal{T} \ra \Func((\mathcal{T}^c)^{\circ},\Vect_k)
\end{gather*}
The previous lemma states that any functor that is cohomological and whose essential image lies in $\vect_k$ is in the essential image of the restricted Yoneda functor. Moreover, we know that the restricted Yoneda functor is full onto the subcategory of cohomological functors taking values in finite-dimensional vector spaces. The next obvious question concerns faithfulness.

\begin{defn}
 A morphism that lies in the kernel of the restricted Yoneda functor is called a \textbf{phantom map}.
\end{defn}

\begin{lem}
 Let $M$ be a bounded above complex of coherent sheaves on $X$. Then, there are no phantom maps from $M$ to a bounded below complex of quasi-coherent sheaves.
\label{lem:phantom}
\end{lem}

\proof We can take a bounded above locally-free coherent resolution $F$ of $M$ and let $F_{\geq i}$ denote the brutal truncation of $F$ at the negative $i$-th step. This is the complex obtained from $F$ by zeroing out the components for $j < -i$. Each $F_{\geq i}$ lies in $D_{\perf}(X)$ and we have a natural map, $\sigma_i: F_{\geq i} \ra F$. Any map, $\phi:M \ra N$, with $N$, a bounded below complex, can be represented by an honest chain map from $F$ to a bounded below complex of injectives $I$. 
\begin{center}
 \leavevmode
 \begin{xy}
  (-60,10)*+{\cdots}="a"; (-45,10)*+{F_{i-1}}="b"; (-30,10)*+{F_i}="c"; (-15,10)*+{F_{i+1}}="d"; (0,10)*+{\cdots}="e"; (15,10)*+{F_{N-1}}="f"; (30,10)*+{0}="g"; (45,10)*+{0}="h"; (60,10)*+{\cdots}="i"; {\ar@{->} "a";"b"}; {\ar@{->} "b";"c"}; {\ar@{->} "c";"d"}; {\ar@{->} "d";"e"}; {\ar@{->} "e";"f"}; {\ar@{->} "f";"g"}; {\ar@{->} "g";"h"}; {\ar@{->} "h";"i"}; (-60,-10)*+{\cdots}="a'"; (-45,-10)*+{0}="b'"; (-30,-10)*+{0}="c'"; (-15,-10)*+{I_{i+1}}="d'"; (0,-10)*+{\cdots}="e'"; (15,-10)*+{I_{N-1}}="f'"; (30,-10)*+{I_N}="g'"; (45,-10)*+{I_{N+1}}="h'"; (60,-10)*+{\cdots}="i'"; {\ar@{->} "a'";"b'"}; {\ar@{->} "b'";"c'"}; {\ar@{->} "c'";"d'"}; {\ar@{->} "d'";"e'"}; {\ar@{->} "e'";"f'"}; {\ar@{->} "f'";"g'"}; {\ar@{->} "g'";"h'"}; {\ar@{->} "h'";"i'"}; {\ar@{->}^0 "b";"b'"}; {\ar@{->}^0 "c";"c'"}; {\ar@{->}^{\phi_{i+1}} "d";"d'"}; {\ar@{->}^{\phi_{N-1}} "f";"f'"}; {\ar@{->}^0 "g";"g'"}; {\ar@{->}^0 "h";"h'"};
 \end{xy}
\end{center}
Since one complex is bounded above and the other is bounded below, $\phi_j$ must be zero outside the interval $[i+1,N-1]$. $\phi \circ \sigma_i$ is null-homotopic if and only if the original map is null-homotopic since any homotopy must also be zero outside $[i+1,N-1]$.
\begin{center}
 \leavevmode
 \begin{xy}
  (-60,10)*+{\cdots}="a"; (-45,10)*+{F_{i-1}}="b"; (-30,10)*+{F_i}="c"; (-15,10)*+{F_{i+1}}="d"; (0,10)*+{\cdots}="e"; (15,10)*+{F_{N-1}}="f"; (30,10)*+{0}="g"; (45,10)*+{0}="h"; (60,10)*+{\cdots}="i"; {\ar@{->} "a";"b"}; {\ar@{->} "b";"c"}; {\ar@{->} "c";"d"}; {\ar@{->} "d";"e"}; {\ar@{->} "e";"f"}; {\ar@{->} "f";"g"}; {\ar@{->} "g";"h"}; {\ar@{->} "h";"i"}; (-60,-10)*+{\cdots}="a'"; (-45,-10)*+{0}="b'"; (-30,-10)*+{0}="c'"; (-15,-10)*+{I_{i+1}}="d'"; (0,-10)*+{\cdots}="e'"; (15,-10)*+{I_{N-1}}="f'"; (30,-10)*+{I_N}="g'"; (45,-10)*+{I_{N+1}}="h'"; (60,-10)*+{\cdots}="i'"; {\ar@{->} "a'";"b'"}; {\ar@{->} "b'";"c'"}; {\ar@{->} "c'";"d'"}; {\ar@{->} "d'";"e'"}; {\ar@{->} "e'";"f'"}; {\ar@{->} "f'";"g'"}; {\ar@{->} "g'";"h'"}; {\ar@{->} "h'";"i'"}; {\ar@{->}^0 "b";"a'"}; {\ar@{->}^0 "c";"b'"}; {\ar@{->}^0 "d";"c'"}; {\ar@{->} "e";"d'"}; {\ar@{->} "f";"e'"}; {\ar@{->}^0 "g";"f'"}; {\ar@{->}^0 "h";"g'"}; {\ar@{->}^0 "i";"h'"};
 \end{xy}
\end{center}
Thus, $\phi \circ \sigma_i: F_{\geq i} \ra N$ is zero in $D(X)$ if and only if $\phi: M \ra N$ is zero in $D(X)$. \qed

The proof says something stronger.

\begin{lem}
 Take a bounded above complex of coherent sheaves, $M$, and a bounded below complex of quasi-coherent sheaves, $N$. There exists a perfect object, $E$, and morphism, $E \ra M$, that induces an isomorphism, $[M,N] \cong [E,N]$, of morphism spaces.
\label{lem:approx}
\end{lem}

The next step is to identify which objects of $D(X)$ give rise to locally-finite functors.

\begin{lem}
 $D^b_{\cohc}(X)$ essentially surjects onto the category of locally-finite cohomological functors via the restricted Yoneda functor.
\label{lem:boundcohcptrep}
\end{lem}

\proof This is essentially theorem A.1 of \cite{BvB02}. First, consider the case of $\P^N_k$. Let $\phi$ be a locally-finite cohomological functor on $D_{\perf}(\P^N_k) \cong D^b_{\coh}(\P^N_k)$. It is represented by a complex $M$. There is an equivalence between $D(\P^N_k)$ and $D(\Mod A)$, which restricts to an equivalence between $D^b_{\coh}(\P^N_k)$ and $D^b(\mod A)$. $A$ is a finite dimensional algebra with finite global dimension. Identify $M$ with its its image. Since $[A,M[j]] \cong H^j(M)$ we see that $M$ has bounded finite-dimensional cohomology and, thus, lies in $D^b(\mod A)$. So $\phi$ is represented by an object of $D^b_{\coh}(\P^N_k)$.

Now consider a locally-finite cohomological functor $\phi$ on $D_{\perf}(X)$. $\phi$ itself is representable by complex, which we will also denote $M$. Choose an embedding $i: X \hookrightarrow \P^N_k$ and consider $\phi' = \phi \circ i^*$.
\begin{gather*}
 \phi'(E) \cong [i^*E,M] \cong [E,i_*M]
\end{gather*}
Thus, $i_*M$ represents $\phi'$ and must lie in $D^b_{\coh}(\P^N_k)$. Since $i: X \ra \P^N_k$ is an embedding,  $i_*: \Qcoh(X) \ra \Qcoh(\P^N_k)$ is exact. In particular, $i_*: D(X) \ra D(\P^N_k)$ is equal to the underived pushforward applied to chain complexes. For a quasi-coherent sheaf, $\mathcal F$, $i_* \mathcal F$ is zero if and only if $\mathcal F$ is zero. Thus, if $i_*M$ has bounded cohomology, the cohomology sheaves of $M$ must be bounded. If $\mathcal F$ is quasi-coherent, then so is $i_*\mathcal F$ since $i_*$ is injective. Thus, $M$ must have coherent cohomology. The support of $i_*M$ is proper as $i_*M$ is a bounded complex of coherent sheaves. Since the support of $M$, viewed as a subset of $\P^N_k$ under $i$, is the support of $i_*M$, the support of $M$ must be proper. Therefore, $M$ lies in $D^b_{\cohc}(X)$. \qed

Lemma \ref{lem:rep} says that the restricted Yoneda functor $D^b_{\coh}(X) \to D_{\perf}(X)^{\vee}$ is full and lemma \ref{lem:phantom} says that it is faithful. The previous lemma says it is essentially surjective. Thus, we have verified the claim of theorem \ref{thm:lfcohom}. \qed

\begin{rmk}
 The conclusion of theorem \ref{thm:lfcohom} also holds if $X$ is Noetherian, possesses the resolution property, and $D(X)$ satisfies an appropriate generation property. Proposition $6.12$ of \cite{Rou03} gives another proof of lemma \ref{lem:boundcohcptrep} using this generation property.
\end{rmk}

Staring at the local finiteness condition, one notices that any perfect object will furnish a locally-finite, homological functor on $D^b_{\cohc}(X)$. Are these all of them? For projective schemes, this question is an easy corollary of results of Rouquier. The following theorem is quite powerful.

\begin{thm}{(Corollary 7.50 of \cite{Rou03}, proposition 5.1.2 of \cite{KMVdB08})}
 Let $X$ be a projective scheme over a field $k$. Any cohomological or homological locally-finite functor on $D^b_{\coh}(X)$ is representable by an object of $D^b_{\coh}(X)$.
\label{thm:Rou}
\end{thm}

\begin{rmk}
 In \cite{Rou03}, Rouquier proves the theorem under the additional assumption that $k$ is perfect. This assumption is removed in \cite{KMVdB08}.
\end{rmk}

\begin{lem}
 An object $A \in D^b_{\coh}(X)$ furnishes a locally-finite functor $[A,-]$ on $D^b_{\cohc}(X)$ if and only if $A$ is perfect.
\label{lem:lffunccoh}
\end{lem}

\proof As noted, if $A$ is perfect, then $[A,-]$ is locally-finite. Assume that $[A,-]$ is locally-finite. Let $x$ be closed point of $X$, $\mathcal{O}_{x,X}$ its local ring, and $\mathcal{O}_x$ the structure sheaf of $x$ in $X$. A bounded complex of finitely-generated $\mathcal{O}_{x,X}$-modules $A$ is quasi-isomorphic to a bounded complex of free modules if and only if 
\begin{displaymath}
 \sum_{i \in \Z} \dim_k [A,\mathcal{O}_x[i]] < \infty.
\end{displaymath}
Any bounded complex of finitely-generated $\mathcal{O}_{x,X}$-modules is quasi-isomorphic to a bounded above complex of finite-rank free $\mathcal{O}_{x,X}$-modules where the differentials are matrices with entries in $m_x$, see \cite{Rob80}. In the case that $A$ is a complex concentrated in degree zero, this is usually called the minimal free resolution of $C$, but the construction works the same for any bounded complex of finitely-generated modules. Applying $[-,\mathcal{O}_x]$ to this minimal complex, we get a complex with zero differentials. If $\sum_{i \in \Z} \dim_k [A,\mathcal{O}_x[i]] < \infty$, the minimal complex must be zero except for finitely many components. Thus, $A$ is quasi-isomorphic to a bounded complex of free $\mathcal{O}_{x,X}$-modules. This quasi-isomorphism extends to some open neighborhood of $x$. Thus, $A$ is perfect. \qed

\begin{prop}
 Let $X$ be a projective scheme over a field $k$. The inclusion, $D_{\perf}(X) \ra \leftexp{\vee}{D^b_{\coh}(X)}$, is an equivalence.
\label{prop:lfhom}
\end{prop}

\proof Theorem \ref{thm:Rou} says that any element of $\leftexp{\vee}{D^b_{\coh}(X)}$ is represntable by an object of $D^b_{\coh}(X)$ and lemma \ref{lem:lffunccoh} says this object must lie in $D_{\perf}(X)$. \qed

\begin{rmk}
 The natural question that now arises is: what happens if $X$ is only quasi-projective? Breaking it down into sub-questions, one wonders: is any locally-finite homological functor, $\phi: D^b_{\cohc}(X) \ra \vect_k$, represented by a perfect complex? Are there any nonzero morphisms between perfect complexes that induce the zero natural transformation as functors on $D^b_{\cohc}(X)$?
\end{rmk}

\section{Pseudo-adjoints}
\label{sec:pseudo-ad}

In this section, unless otherwise stated, $X$ and $Y$ will be projective schemes over a field $k$.

It is important to keep in mind, when contemplating the results of the previous section, that there are two components of representability. The first is representing a functor by an object of an appropriate category. The second is representing natural transformations as morphisms between the representing objects. In the case that the appropriate category is actually the underlying category, the Yoneda lemma makes quick work of the second issue. However, in the cases of interest in this paper, the appropriate category is not the underlying category. Thus, the second issue must be properly addressed.

It is the representability of natural transformations that is the engine behind most results related to representability in the usual setting. For instance, general theorems giving sufficient conditions for representability of a functor have easy corollaries involving the existence of adjoints. In this section, we will collect some easy corollaries, many involving (pseudo-)adjunction, which follow from the representability results in the previous section.

\begin{defn}
 If $F: D_{\perf}(X) \ra D_{\perf}(Y)$ is an exact functor, a \textbf{right pseudo-adjoint} to $F$ is a functor $F^{\vee}: D^b_{\coh}(Y) \ra D^b_{\coh}(X)$ so that we have natural isomorphisms
\begin{displaymath}
 [F(A),B] \cong [A,F^{\vee}(B)]
\end{displaymath}
for any pair of objects $A$ in $D_{\perf}(X)$ and $B$ in $D^b_{\coh}(Y)$. 

If $G: D^b_{\coh}(X) \ra D^b_{\coh}(Y)$ is an exact functor, a \textbf{left pseudo-adjoint} to $G$ is a functor $\leftexp{\vee}{G}: D_{\perf}(Y) \ra D_{\perf}(X)$ so that we have natural isomorphisms
\begin{displaymath}
 [A,G(B)] \cong [\leftexp{\vee}{G}(A),B]
\end{displaymath}
for any pair of objects $A$ in $D_{\perf}(Y)$ and $B$ in $D^b_{\coh}(X)$.
\end{defn}

\begin{rmk}
 We can also extend the definition appropriately to the case where $X$ and $Y$ are quasi-projective. Hopefully, it will be clear which results in this section extend, after appropriate modification, to quasi-projective schemes over $k$.
\end{rmk}

\begin{prop}
 Any exact functor $F: D_{\perf}(X) \ra D_{\perf}(Y)$ possesses a right pseudo-adjoint, which is unique. Any exact functor $G: D^b_{\coh}(X) \ra D^b_{\coh}(Y)$ possesses a left pseudo-adjoint, which is unique.
\label{prop:exisuniq}
\end{prop}

\proof The proofs for both statements are similar. By precomposition, $F$ induces a functor $D_{\perf}(Y)^{\vee} \ra D_{\perf}(X)^{\vee}$. Under the equivalence in theorem \ref{thm:lfcohom}, this uniquely specifies $F^{\vee}$. The statement for left pseudo-adjoints follows via a completely analogously argument with the use of proposition \ref{prop:lfhom} in place of theorem \ref{thm:lfcohom}. \qed

The following ``double dualization'' statement is an immediate consequence of uniqueness of pseudo-adjoints.

\begin{cor}
 We have natural isomorphisms
\begin{align*}
 \leftexp{\vee}{\left(F^{\vee}\right)} & \cong F \\
 \left(\leftexp{\vee}{G}\right)^{\vee} & \cong G.
\end{align*}
\end{cor}

\begin{lem}
 $F: D_{\perf}(X) \ra D_{\perf}(Y)$ possesses a left adjoint if and only if $F$ admits an extension to $D^b_{\coh}(X)$ which takes values in $D^b_{\coh}(Y)$. If $G: D^b_{\coh}(X) \ra D^b_{\coh}(Y)$ posseses a right adjoint, $G$ takes perfect complexes to perfect complexes.
\label{lem:adj->padj}
\end{lem}

\proof If $F: D_{\perf}(X) \ra D_{\perf}(Y)$ possesses a left adjoint, $F^*$, then $[A,F(B)] \cong [F^*(A),B] \cong [A,(F^*)^{\vee}(B)]$ for $A \in D_{\perf}(Y)$ and $B \in D_{\perf}(X)$. Since $(F^*)^{\vee}(B)$ and $F(B)$ furnish the same cohomological functor on $D_{\perf}(Y)$, they must be isomorphic. This isomorphism is natural in $B$. So $(F^*)^{\vee}$ is an extension of $F$ to $D^b_{\coh}(X)$ which takes values in $D^b_{\coh}(Y)$. Now suppose that $F$ admits such an extension, call it $\tilde{F}$. Then, for $A \in D_{\perf}(Y)$ and $B \in D_{\perf}(B)$, we have $[A,F(B)] \cong [\leftexp{\vee}{(\tilde{F})}(A),B]$ by pseudo-adjunction. So we get an honest adjunction.

If $G$ has a right adjoint, $G^!$, then, as $[G(A),B] \cong [A,G^!(B)]$ for $A \in D_{\perf}(X)$ and $B \in D^b_{\coh}(Y)$, $[G(A),-]$ is a locally-finite cohomological functor on $D^b_{\coh}(Y)$. Consequently, $G(A)$ is perfect.
\qed

\begin{rmk}
 To prove that $G$ possesses a right adjoint if it takes perfect complexes to perfect complexes would require us to extend the natural isomorphisms provided by pseudo-adjunction. One can represent $[C,-]$, for a bounded complex of coherent sheaves, as the limit of maps out of truncations of $C$. Using this, one can extend the natural transformations. We leave the details to the reader.
\end{rmk}

The following allows us to clarify the relationship between autoequivalences of perfect derived categories and autoequivalences of bounded derived categories.  

\begin{lem}
 If $F$ is an equivalence, then so is $F^{\vee}$. If $G$ is an equivalence, then so is $\leftexp{\vee}{G}$. Consequently, we see that two projective schemes have equivalent perfect derived categories if and only if they have equivalent bounded derived categories. Moreover, for a projective scheme, $X$, the groups of autoequivalences of $D_{\perf}(X)$ and the group of auto-equivalences of $D^b_{\coh}(X)$ are isomorphic. 
\label{lem:misclfdualresults}
\end{lem}

\proof From proof of proposition \ref{prop:exisuniq}, we know that $F^{\vee}$ is isomorphic to the functor $D_{\perf}(Y)^{\vee} \to D_{\perf}(X)^{\vee}$ given by precomposition by $F$. This is under the equivalences $D_{\perf}(X)^{\vee} \cong D^b_{\coh}(X)$ and $D_{\perf}(Y)^{\vee} \cong D^b_{\coh}(Y)$ from theorem \ref{thm:lfcohom}. From this description, it is clear that $F^{\vee}$ is an equivalence if $F$ is. The argument for $\leftexp{\vee}{G}$ is completely analogous. The final statements follow immediately from the first two. \qed

%
%


One can restrict an auto-equivalence of $D(X)$ to get an auto-equivalence of $D_{\perf}(X)$. There is an inverse to this restriction, see \cite{Bal09b} or \cite{LO09}.

We have the following notion due to Buchweitz, \cite{Buc86}, for affine varieties and Orlov, \cite{Orl04}, in general.

\begin{defn}
 Let $X$ be a projective scheme over a field. The Verdier quotient
\begin{displaymath}
 D^b_{\coh}(X)/D_{\perf}(X)
\end{displaymath}
 is called the \textbf{triangulated category of singularities} of $X$ and is denoted by $D_{\sing}(X)$.
\end{defn}

\begin{cor}
 If $X$ and $Y$ are projective schemes over a field with equivalent bounded derived categories of coherent sheaves or perfect derived categories, then they have equivalent categories of singularities. In particular, $X$ is regular if and only if $Y$ is regular.
\end{cor}

\proof Lemma \ref{lem:misclfdualresults} tells us that the perfect derived categories are equivalent if and only if the bounded derived categories are equivalent. Moreover, the equivalence of the bounded derived categories must restrict to an equivalence of the perfect derived categories. Therefore, this equivalence induces an equivalence between the quotient categories. It is easy corollary of Serre's criteria for regularity of local ring that the singularity category is zero if and only if the scheme is regular. \qed

\begin{rmk}
 The regularity statement follows immediately from the main result of \cite{Orl04}. The statement that equivalent bounded derived categories of coherent sheaves implies equivalent categories of singularities is implicit in corollary 1.12 of \cite{Orl06}. 
\end{rmk}

The triangulated structure on $D^b_{\coh}(X)$ (and $D^b_{\coh}(Y)$) is not so easy to capture in terms of locally-finite cohomological functors. Indeed, this is the most common manifestation of the non-functoriality  of cones in triangulated categories. Nevertheless, we have the following result.

\begin{lem}
 If $F: D_{\perf}(X) \ra D_{\perf}(Y)$ is an exact functor, then $F^{\vee}$ is also exact.
\label{lem:pseudoadjointexact}
\end{lem}

\proof We can define the shift functor on $D^b_{\coh}(X)$ (and $D^b_{\coh}(Y)$) purely in terms of locally-finite cohomological functors as the functor that sends $\Phi$ to $\Phi \circ [-1]$. $F^{\vee}$ is precomposition with $F$ and thus commutes with shifts since $F$ is exact.

Take an exact triangle
\begin{center}
 \leavevmode
 \begin{xy}
  (-10,10)*+{A}="a"; (10,10)*+{B}="b"; (0,-5)*+{C}="c"; {\ar@{->}^{\alpha} "a";"b"}; {\ar@{->}^{\beta} "b";"c"}; {\ar@{-->}^{[1]} "c";"a"}
 \end{xy}
\end{center}
in $D^b_{\coh}(Y)$. Let $Z$ be a cone over $F^{\vee}(\alpha)$. We seek a morphism $\lambda: Z \ra F^{\vee}(C)$ making
\begin{center}
 \leavevmode
 \begin{xy}
  (-45,10)*+{F^{\vee}(A)}="a"; (-15,10)*+{F^{\vee}(B)}="b"; (15,10)*+{Z}="c"; (45,10)*+{F^{\vee}(A)[1]}="d"; {\ar@{->}^{F^{\vee}(\alpha)} "a";"b"}; {\ar@{->} "b";"c"}; {\ar@{->} "c";"d"};
  (-45,-10)*+{F^{\vee}(A)}="a'"; (-15,-10)*+{F^{\vee}(B)}="b'"; (15,-10)*+{F^{\vee}(C)}="c'"; (45,-10)*+{F^{\vee}(A)[1]}="d'"; {\ar@{->}^{F^{\vee}(\alpha)} "a'";"b'"}; {\ar@{->} "b'";"c'"}; {\ar@{->} "c'";"d'"};
  {\ar@{->}^= "a";"a'"}; {\ar@{->}^= "b";"b'"}; {\ar@{->}^{\lambda} "c";"c'"}; {\ar@{->}^= "d";"d'"};
 \end{xy}
\end{center}
commute where the bottom row results from application of $F^{\vee}$ to the triangle $A \to B \to C$. 

Let us first work assuming we know such a $\lambda$ exists. Let us prove the following lemma.

\begin{lem}
 Let $\mathcal T$ be a triangulated category. Assume we have a commutative diagram of the form
\begin{center}
 \leavevmode
 \begin{xy}
  (-45,10)*+{X}="a"; (-15,10)*+{Y}="b"; (15,10)*+{Z}="c"; (45,10)*+{X[1]}="d"; {\ar@{->} "a";"b"}; {\ar@{->} "b";"c"}; {\ar@{->} "c";"d"};
  (-45,-10)*+{X}="a'"; (-15,-10)*+{Y}="b'"; (15,-10)*+{W}="c'"; (45,-10)*+{X[1]}="d'"; {\ar@{->} "a'";"b'"}; {\ar@{->} "b'";"c'"}; {\ar@{->} "c'";"d'"};
  {\ar@{->}^{\cong} "a";"a'"}; {\ar@{->}^{\cong} "b";"b'"}; {\ar@{->}^{\lambda} "c";"c'"}; {\ar@{->}^{\cong} "d";"d'"};
 \end{xy}
\end{center}
 where the vertical maps, except for possibly $\lambda$, are isomorphisms. If $E$ is an object of $\mathcal T$ for which application of $[E,-]$, respectively $[-,E]$, to the rows yields exact sequences, then $[E,\lambda]$, respectively $[\lambda,E]$, is an isomorphism. If $[W,-]$ and $[Z,-]$ or $[-,W]$ and $[-,Z]$ convert the rows into exact sequences, then $W$ is isomorphic to $Z$.
\label{lem:diagram}
\end{lem}

\proof Let $E$ be an object for which the rows in the diagram
\begin{center}
 \leavevmode
 \begin{xy}
  (-45,10)*+{[E,X]}="a"; (-15,10)*+{[E,Y]}="b"; (15,10)*+{[E,Z]}="c"; (45,10)*+{[E,X[1]]}="d"; {\ar@{->} "a";"b"}; {\ar@{->} "b";"c"}; {\ar@{->} "c";"d"};
  (-45,-10)*+{[E,X]}="a'"; (-15,-10)*+{[E,Y]}="b'"; (15,-10)*+{[E,W]}="c'"; (45,-10)*+{[E,X[1]]}="d'"; {\ar@{->} "a'";"b'"}; {\ar@{->} "b'";"c'"}; {\ar@{->} "c'";"d'"};
  {\ar@{->}^{\cong} "a";"a'"}; {\ar@{->}^{\cong} "b";"b'"}; {\ar@{->}^{[E,\lambda]} "c";"c'"}; {\ar@{->}^{\cong} "d";"d'"};
 \end{xy}
\end{center}
 are exact. The $5$-lemma implies that $[E,\lambda]$ is an isomorphism. If $[-,E]$ converts the rows into exact sequences, then the argument is similar. Let us assume that $Z$ and $W$ convert the rows into exact sequences. There is an element of $[Z,W]$ whose composition with $\lambda$ is $\operatorname{id}_Z$ and an element of $[W,Z]$ whose composition with $\lambda$ is $\operatorname{id}_W$. Consequently, $\lambda$ is an isomorphism. If $[-,W]$ and $[-,Z]$ convert the rows into exact sequences, then the argument is similar. \qed

To show $\lambda$ is an isomorphism, it is sufficient, by theorem \ref{thm:lfcohom}, to show that it induces an isomorphism after application of $[E,-]$ for each perfect $E$. Let $E$ be perfect. Applying $[E,-]$ to the diagram
\begin{center}
 \leavevmode
 \begin{xy}
  (-45,10)*+{F^{\vee}(A)}="a"; (-15,10)*+{F^{\vee}(B)}="b"; (15,10)*+{Z}="c"; (45,10)*+{F^{\vee}(A)[1]}="d"; {\ar@{->}^{F^{\vee}(\alpha)} "a";"b"}; {\ar@{->} "b";"c"}; {\ar@{->} "c";"d"};
  (-45,-10)*+{F^{\vee}(A)}="a'"; (-15,-10)*+{F^{\vee}(B)}="b'"; (15,-10)*+{F^{\vee}(C)}="c'"; (45,-10)*+{F^{\vee}(A)[1]}="d'"; {\ar@{->}^{F^{\vee}(\alpha)} "a'";"b'"}; {\ar@{->} "b'";"c'"}; {\ar@{->} "c'";"d'"};
  {\ar@{->}^= "a";"a'"}; {\ar@{->}^= "b";"b'"}; {\ar@{->}^{\lambda} "c";"c'"}; {\ar@{->}^= "d";"d'"};
 \end{xy}
\end{center}
converts the top row to an exact sequence as $F^{\vee}(A) \to F^{\vee}(B) \to Z$ is a distinguished triangle. Via pseudo-adjunction, the bottom row is isomorphic to
\begin{displaymath}
 [F(E),A] \to [F(E),B] \to [F(E),C] \to [F(E),A[1]].
\end{displaymath}
Since this sequence is exact, the bottom row is also exact. Appealing to lemma \ref{lem:diagram} shows that $[E,\lambda]$ is an isomorphism whenever $E$ is perfect.

To establish the existence of the desired $\lambda$, we want to mimic the early steps of the proof of exactness of adjoints in proposition 1.4 of \cite{BK89}. However, we cannot use pseudo-adjunction unless the domain of a morphism is perfect. To deal with this, we approximate complexes by perfect ones.

Let $E_i^{F^{\vee}(A)}$ denote the brutal truncation of a locally-free resolution of $F^{\vee}(A)$ at the $-i$th stage and let $E_j^{F^{\vee}(B)}$ denote the brutal truncation of a locally-free resolution of $F^{\vee}(B)$ at the $-j$th stage. $E_i^{F^{\vee}(A)}$ is compact and $F^{\vee}(B)$ is the homotopy colimit of the $E_j^{F^{\vee}(B)}$. By Lemma \ref{lem:hocolim/colim}, for large $i$, there exists a map $[F^{\vee}(\alpha)]_i$ making
\begin{center}
 \leavevmode
 \begin{xy}
  (-15,10)*+{E_i^{F^{\vee}(A)}}="a"; (15,10)*+{E_{j(i)}^{F^{\vee}(B)}}="b"; {\ar@{->}^{[F^{\vee}(\alpha)]_i} "a";"b"};
  (-15,-10)*+{F^{\vee}(A)}="a'"; (15,-10)*+{F^{\vee}(B)}="b'"; {\ar@{->}^{F^{\vee}(\alpha)} "a'";"b'"}; {\ar@{->} "a";"a'"}; {\ar@{->} "b";"b'"};
 \end{xy}
\end{center}
commute. Note that we are always free to increase the value of $j(i)$ without affecting the result.
Complete this commutative square to a morphism of triangles
\begin{center}
 \leavevmode
 \begin{xy}
  (-45,10)*+{E_i^{F^{\vee}(A)}}="a"; (-15,10)*+{E_{j(i)}^{F^{\vee}(B)}}="b"; (15,10)*+{C_{i}}="c"; (45,10)*+{E_i^{F^{\vee}(A)}[1]}="d"; {\ar@{->}^{[F^{\vee}(\alpha)]_{i}} "a";"b"}; {\ar@{->} "b";"c"}; {\ar@{->} "c";"d"}; 
  (-45,-10)*+{F^{\vee}(A)}="a'"; (-15,-10)*+{F^{\vee}(B)}="b'"; (15,-10)*+{Z}="c'"; (45,-10)*+{F^{\vee}(A)[1]}="d'"; {\ar@{->}^{F^{\vee}(\alpha)} "a'";"b'"}; {\ar@{->} "b'";"c'"}; {\ar@{->} "c'";"d'"};
  {\ar@{->} "a";"a'"}; {\ar@{->} "b";"b'"}; {\ar@{->} "d";"d'"}; {\ar@{->} "c";"c'"};
 \end{xy}
\end{center}
By naturality of pseudo-adjunction, we have a commutative square
\begin{center}
 \leavevmode
 \begin{xy}
  (-15,10)*+{F(E_i^{F^{\vee}(A)})}="a"; (15,10)*+{F(E_{j(i)}^{F^{\vee}(B)})}="b"; {\ar@{->} "a";"b"};
  (-15,-10)*+{A}="a'"; (15,-10)*+{B}="b'"; {\ar@{->} "a'";"b'"}; {\ar@{->} "a";"a'"}; {\ar@{->} "b";"b'"};
 \end{xy}
\end{center}
We can complete this to a morphism of triangles
\begin{center}
 \leavevmode
 \begin{xy}
  (-45,10)*+{F(E_i^{F^{\vee}(A)})}="a"; (-15,10)*+{F(E_{j(i)}^{F^{\vee}(B)})}="b"; (15,10)*+{F(C_{i})}="c"; (45,10)*+{F(E_i^{F^{\vee}(A)})[1]}="d"; {\ar@{->} "a";"b"}; {\ar@{->} "b";"c"}; {\ar@{->} "c";"d"}; 
  (-45,-10)*+{A}="a'"; (-15,-10)*+{B}="b'"; (15,-10)*+{C}="c'"; (45,-10)*+{A[1]}="d'"; {\ar@{->} "a'";"b'"}; {\ar@{->} "b'";"c'"}; {\ar@{->} "c'";"d'"}; {\ar@{->} "a";"a'"}; {\ar@{->} "b";"b'"}; {\ar@{->} "d";"d'"}; {\ar@{->} "c";"c'"};
 \end{xy}
\end{center}
By pseudo-adjunction, we get a map $\theta: C_i \ra F^{\vee}(C)$. We now have a slightly larger commutative diagram.
\begin{center}
 \leavevmode
 \begin{xy}
  (-45,15)*+{E_i^{F^{\vee}(A)}}="a"; (-15,15)*+{E_{j(i)}^{F^{\vee}(B)}}="b"; (15,15)*+{C_{i}}="c"; (45,15)*+{E_i^{F^{\vee}(A)}[1]}="d"; {\ar@{->}^{[F^{\vee}(\alpha)]_{i}} "a";"b"}; {\ar@{->} "b";"c"}; {\ar@{->} "c";"d"}; 
  (-45,0)*+{F^{\vee}(A)}="a'"; (-15,0)*+{F^{\vee}(B)}="b'"; (15,0)*+{Z}="c'"; (45,0)*+{F^{\vee}(A)[1]}="d'";
  (-45,-15)*+{F^{\vee}(A)}="a''"; (-15,-15)*+{F^{\vee}(B)}="b''"; (15,-15)*+{F^{\vee}(C)}="c''"; (45,-15)*+{F^{\vee}(A)[1]}="d''"; {\ar@{->}^{F^{\vee}(\alpha)} "a'";"b'"}; {\ar@{->} "b'";"c'"}; {\ar@{->} "c'";"d'"};
  {\ar@{->}^{F^{\vee}(\alpha)} "a''";"b''"}; {\ar@{->} "b''";"c''"}; {\ar@{->} "c''";"d''"};
  {\ar@{->} "a";"a'"}; {\ar@{->} "b";"b'"}; {\ar@{->} "d";"d'"}; {\ar@{->} "c";"c'"};
  {\ar@{->}^{=} "a'";"a''"}; {\ar@{->}^{=} "b'";"b''"}; {\ar@{->}^{=} "d'";"d''"}; {\ar@/^1pc/ "c";"c''"};
  (20.5,3)*+{\theta};
 \end{xy}
\end{center}
The top squares commute because we have a morphism of triangles. The bottom square commutes tautologically and any square involving $\theta$ commutes by naturality of pseudo-adjunction. We just need to fill in the missing map from $Z$ to $F^{\vee}(C)$.

For any bounded complex of coherent sheaves, $V$, make a choice of locally-free resolution of $V$, whose brutal truncation at the $-i$th step will be denoted by $E_i^V$. By lemma \ref{lem:approx}, if we are given another bounded complex of coherent sheaves, $W$, there exists an $N_0$ such that for $i > N_0$ the truncation map, $E_i^V \ra V$, induces an isomorphism
\begin{displaymath}
 [V,W] \overset{\sim}{\longrightarrow} [E_i^V,W]
\end{displaymath}
of the morphism vector spaces. 
Choose $W = F^{\vee}(C)$ and $i$ large enough so that restrictions induces isomorphisms $[E_i^{F^{\vee}(A)},F^{\vee}(C)] \cong [F^{\vee}(A),F^{\vee}(C)]$ and $[E_{j(i)}^{F^{\vee}(B)},F^{\vee}(C)] \cong [F^{\vee}(B),F^{\vee}(C)]$. By the $5$-lemma applied to the upper part of the diagram, the map $C_i \ra Z$ induces an isomorphism $[C_i,F^{\vee}(C)] \cong [Z,F^{\vee}(C)]$. Therefore, $\theta$ is the restriction of unique map $\lambda: Z \ra F^{\vee}(C)$. To see that $\lambda: Z \ra F^{\vee}(C)$ makes the total diagram commute, let us label the maps $\xi: 
E_{j(i)}^{F^{\vee}(B)} \ra C_i$, $r_{j(i)}: E_{j(i)}^{F^{\vee}(B)} \ra F^{\vee}(B)$, and $\nu:F^{\vee}(B) \ra Z$. Equality of $\lambda \circ \nu$ and $F^{\vee}(\beta)$ is now equivalent to equality of their restrictions to the top row. But, their restrictions are $\theta \circ \xi$ and $F^{\vee}(\beta) \circ r_{j(i)}$, respectively. We have already seen these are equal. To show that the other square involving $\lambda$ commutes, we can increase $i$ so that, additionally, restrictions induce isomorphisms $[E_i^{F^{\vee}(A)},F^{\vee}(A)[1]] \cong [F^{\vee}(A),F^{\vee}(A)[1]]$ and $[E_{j(i)}^{F^{\vee}(B)},F^{\vee}(A)[1]] \cong [F^{\vee}(B),F^{\vee}(A)[1]]$ and repeat the argument. Thus, we have found the desired map $\lambda: Z \ra F^{\vee}(C)$. \qed

The following lemma is a little easier.

\begin{lem}
 If $F:D^b_{\coh}(X) \ra D^b_{\coh}(Y)$ is an exact functor, then its left pseudo-adjoint, $\leftexp{\vee}{F}: D_{\perf}(Y) \ra D_{\perf}(X)$, is also exact.
\end{lem}

\proof Take an exact triangle
\begin{center}
 \leavevmode
 \begin{xy}
  (-10,10)*+{A}="a"; (10,10)*+{B}="b"; (0,-5)*+{C}="c"; {\ar@{->}^{\alpha} "a";"b"}; {\ar@{->}^{\beta} "b";"c"}; {\ar@{-->}^{[1]} "c";"a"}
 \end{xy}
\end{center}
in $D^b_{\perf}(Y)$. Let $Z$ be a cone over $\leftexp{\vee}{F}(\alpha)$. We seek a morphism $\lambda:\leftexp{\vee}{F}(C) \to Z$ making
\begin{center}
 \leavevmode
 \begin{xy}
  (-45,10)*+{\leftexp{\vee}{F}(A)}="a"; (-15,10)*+{\leftexp{\vee}{F}(B)}="b"; (15,10)*+{\leftexp{\vee}{F}(C)}="c"; (45,10)*+{\leftexp{\vee}{F}(A)[1]}="d"; {\ar@{->}^{\leftexp{\vee}{F}(\alpha)} "a";"b"}; {\ar@{->} "b";"c"}; {\ar@{->} "c";"d"};
  (-45,-10)*+{\leftexp{\vee}{F}(A)}="a'"; (-15,-10)*+{\leftexp{\vee}{F}(B)}="b'"; (15,-10)*+{Z}="c'"; (45,-10)*+{\leftexp{\vee}{F}(A)[1]}="d'"; {\ar@{->}^{\leftexp{\vee}{F}(\alpha)} "a'";"b'"}; {\ar@{->} "b'";"c'"}; {\ar@{->} "c'";"d'"};
  {\ar@{->}^= "a";"a'"}; {\ar@{->}^= "b";"b'"}; {\ar@{->}^{\lambda} "c";"c'"}; {\ar@{->}^= "d";"d'"};
 \end{xy}
\end{center}
commute where the top row results from application of $\leftexp{\vee}{F}$ to the triangle $A \to B \to C$. Given such a $\lambda$, we can apply lemma \ref{lem:diagram} to conclude that $\lambda$ is an isomorphism. 

Let us now show that such a $\lambda$ exists. Apply $F$ to the bottom row and note, for $D \in D_{\perf}(Y)$, we have canonical maps, $\gamma_D$, in $[D,F \circ \leftexp{\vee}{F}(D)] \cong [\leftexp{\vee}{F}(D),\leftexp{\vee}{F}(D)]$ corresponding to $\operatorname{id}_{\leftexp{\vee}{F}(D)}$. We can find a morphism $\gamma: C \to F(Z)$ making the diagram
\begin{center}
 \leavevmode
 \begin{xy}
  (-45,10)*+{A}="a"; (-15,10)*+{B}="b"; (15,10)*+{C}="c"; (45,10)*+{A[1]}="d"; {\ar@{->}^{\alpha} "a";"b"}; {\ar@{->} "b";"c"}; {\ar@{->} "c";"d"};
  (-45,-10)*+{F\circ\leftexp{\vee}{F}(A)}="a'"; (-15,-10)*+{F\circ\leftexp{\vee}{F}(B)}="b'"; (15,-10)*+{FZ}="c'"; (45,-10)*+{F\circ\leftexp{\vee}{F}(A)[1]}="d'"; {\ar@{->}^{F \circ F^{\vee}(\alpha)} "a'";"b'"}; {\ar@{->} "b'";"c'"}; {\ar@{->} "c'";"d'"};
  {\ar@{->}^{\gamma_A} "a";"a'"}; {\ar@{->}^{\gamma_B} "b";"b'"}; {\ar@{->}^{\gamma} "c";"c'"}; {\ar@{->}^{\gamma_{A[1]}} "d";"d'"};
 \end{xy}
\end{center}
commute. By naturality of pseudo-adjunction, $\gamma: C \to F(Z)$ corresponds to the desired $\lambda: \leftexp{\vee}{F}(C) \to Z$.

\section{Rouquier functors}
\label{sec:Roufunc}
 
In section \ref{sec:recon}, we wish to demonstrate the utility of the ideas and results found in the previous sections by extending the reconstruction result of Bondal and Orlov to Gorenstein projective varieties. Before we can do this, we need one more idea, that of a Rouquier functor.

Let us recall the definition of a Serre functor.
\begin{defn}
 Let $\mathcal{C}$ be a $k$-linear category. A \textbf{weak Serre functor} $S: \mathcal{C} \ra \mathcal{C}$ is an endofunctor for which there are natural isomorphisms
\begin{displaymath}
 \eta_{A,B}:[B,S(A)] \ra [A,B]^*
\end{displaymath}
for any pair of objects $A$ and $B$ of $\mathcal{C}$. $S$ is a \textbf{Serre functor} if, in addition, it is an autoequivalence.
\end{defn}

\begin{rmk}
 This definition is slightly different than the one that commonly appears in the literature, see \cite{BK89}. We have reversed the dualization. Of course, if the category has finite dimensional morphism spaces, the two definitions are equivalent.
\end{rmk}

\begin{eg}
 The canonical example for a Serre functor is the following: let $X$ be a smooth projective variety over a field $k$ and let $\omega_X$ be the canonical bundle. Then, the Serre functor on $D^b_{\coh}(X)$ is $-\otimes \omega_X[\dim X]$. This follows from Serre duality.
\end{eg}

\begin{defn}
 Let $\mathcal{C}$ and $\mathcal{D}$ be $k$-linear categories and $F: \mathcal{C} \ra \mathcal{D}$ a $k$-linear functor. A \textbf{Rouquier functor} for $F$ is a $k$-linear functor, $R_F: \mathcal{C} \ra \mathcal{D}$, for which there are natural isomorphisms
\begin{displaymath}
 \eta_{A,B}: [B,R_F(A)] \ra [F(A),B]^*
\end{displaymath}
\end{defn}

\begin{rmk}
 Recall that naturality of $\eta$ is equivalent to the following statements: let $A,A'$ be objects of $\mathcal C$ and $B,B'$ be objects of $\mathcal D$. Let $\lambda \in [B,R_F(A)]$, $\rho \in [F(A'),B]$, and $\tau \in [F(A),B']$. Let $\alpha \in [A,A']$ and $\beta \in [B',B]$. We have equalities
\begin{align*}
 \eta_{A',B}(R_F(\alpha) \circ \lambda)(\rho) & = \eta_{A,B}(\lambda)(\rho \circ F(\alpha)) \\
 \eta_{A,B'}(\lambda \circ \beta)(\tau) &= \eta_{A,B}(\lambda)(\beta \circ \tau).
\end{align*}
 We will use these equalities repeatedly in the proof of proposition \ref{prop:Roufunctriang}.
\end{rmk}


\begin{rmk}
 The definition of a Rouquier functor should be viewed as an attempt to relativize the notion of Serre functor. Indeed, if we take $F$ to be the identity functor, $R_F$ is just a weak Serre functor.
\end{rmk}

The definition is natural so we should expect most categorical notions to respect Rouquier functors.

\begin{lem}
 A necessary and sufficient condition for the existence of $R_F$ is the representablity of the functors $[F(A),-]^*$ for objects $A$ of $\mathcal{C}$. If $R_F$ exists, it is unique.
\label{lem:condforRouexist}
\end{lem}

\proof Clearly, if $R_F$ exists, then $[F(A),-]^*$ is representable. Indeed, it is represented by $R_F(A)$. Assume that $[F(A),-]^*$ is representable for all $A$. We set $R_F(A)$ to be a choice of the representing object for $[F(A),-]^*$. If we have a morphism $\alpha: A \ra A'$, this gives a natural transformation of functors $F(\alpha):[F(A),-]^* \ra [F(A'),-]^*$ which corresponds to a natural transformation $[-,R_F(A)]$ $\ra$ $[-,R_F(A')]$ and, by the Yoneda lemma, a morphism $R_F(\alpha): R_F(A) \ra R_F(A')$. 

Assume we have two choices of Rouquier functors, $R_F$ and $R_F'$. As $[-,R_F(A)]$ and $[-,R_F'(A)]$ are both isomorphic to $[F(A),-]^*$, we see that there is a isomorphism $\nu_A: R_F(A) \to R_F'(A)$ which is natural in $A$. The resulting natural transformation $\nu: R_F \to R_F'$ is an isomorphism.
\qed

The author would like to thank the referee for pointing out the following useful fact.

\begin{lem}
 Let $G: \mathcal D \to \mathcal E$ be an exact $k$-linear functors between $k$-linear triangulated categories. If $R_G$ exists, then $R_{G \circ F}$ exists and is isomorphic to $R_G \circ F$. If $G$ is an equivalence and either $R_F$ or $R_{G \circ F}$ exists, then $R_{G \circ F} \cong G \circ R_F$.
\label{lem:Rouqcompose}
\end{lem}

\proof If $R_G$ exists, we have the following natural isomorphism
\begin{displaymath}
 [G \circ F(A),B]^* \cong [B,R_G \circ F(A)]
\end{displaymath}
 for any pair of objects $A \in \mathcal C$ and $B \in \mathcal E$. By lemma \ref{lem:condforRouexist}, $R_{G \circ F}$ exists and is isomorphic to $R_G \circ F$.

 If $G$ is an equivalence and $R_F$ exists, we have the following natural isomorphisms
\begin{displaymath}
 [G \circ F(A),B]^* \cong [F(A),G^{-1}(B)]^* \cong [G^{-1}(B),R_F(A)] \cong [B,G \circ R_F(A)]. 
\end{displaymath}
 By lemma \ref{lem:condforRouexist}, $R_{G \circ F}$ exists and is isomorphic to $G \circ R_F$. If $R_{F \circ G}$ exists, the above shows that $[F(A),-]^*$ is representable for each $A$ so by lemma \ref{lem:condforRouexist} $R_F$ exists and we have $R_{G \circ F} \cong G \circ R_F$. \qed

\begin{lem}
 If $\mathcal{D}$ has finite-dimensional morphism spaces and $F$ is fully-faithful, then $R_F$ is fully-faithful.
\label{lem:Roufullfaith}
\end{lem}

\proof We have natural isomorphisms 
\begin{displaymath}
 [R_F(A),R_F(B)] \cong [F(B),R_F(A)]^* \cong [F(A),F(B)]^{**} \cong [F(A),F(B)] \cong [A,B].
\end{displaymath}
\qed

There is one useful case where we can guarantee the existence of a Rouquier functor. The following result and lemma \ref{lem:RouisSerre} were first observed in \cite{Rou03}, which inspired the name, Rouquier functor.

\begin{lem}
 Let $\mathcal{T}$ be a compactly-generated triangulated category. Then, the Rouquier functor for the inclusion $\mathcal{T}^c \hookrightarrow \mathcal{T}$ exists.
\end{lem}

\proof By lemma \ref{lem:condforRouexist}, we just need to show that $[A,-]^*$ is representable for each compact $A$. $[A,-]^*$ takes coproducts to products and is, thus, representable by Brown representability. \qed \\
We denote the Rouquier functor in this case by $R_{\mathcal{T}}$.

As we know that the Rouquier functor for $\mathcal T^c \hookrightarrow \mathcal T$ always exists if $\mathcal T$ is compactly-generated, lemma \ref{lem:Rouqcompose} gives us the following.

\begin{lem}
 Let $F$ be a functor from $\mathcal{S}$ to $\mathcal{T}$. If $\mathcal{T}$ is compactly generated and the essential image of $F$ lies in $\mathcal{T}^c$, then $R_F$ exists.
\end{lem}

Again, we return to the geometric setting for the main example.

\begin{eg}
 Let $X$ be a quasi-projective scheme over a field $k$. Let $f: X \ra \Spec k$ be the structure map and $f^!$ the right adjoint to $f_*$. The Rouquier functor for the inclusion, $D_{\perf}(X) \hookrightarrow D(X)$, is $- \otimes f^!\mathcal{O}_{\Spec k}$. Let us denote this functor by $R_X$. To verify that $R_X$ is indeed the Rouquier functor, take $A$ to be a perfect complex and $B$ to be a complex of quasi-coherent sheaves on $X$ and observe the following natural isomorphisms:
\begin{gather*}
 [A,B]^* \cong [\mathcal{O}_X,A^{\vee} \otimes B]^* = [f_*(A^{\vee} \otimes B),\mathcal{O}_k] \\ \cong [A^{\vee} \otimes B,f^!\mathcal{O}_{\Spec k}] \cong [B,A \otimes f^!\mathcal{O}_{\Spec k}]
\end{gather*}
 If $X$ is projective, then $[A,-]^*$ is a locally-finite cohomological functor. So $R_X$ maps $D_{\perf}(X)$ into $D^b_{\coh}(X)$.
\label{eg:rou}
\end{eg}

\begin{lem}
 Let $\mathcal{T}$ be a compactly-generated triangulated category. If $\mathcal{T}^c$ has a weak Serre functor $S$, it is isomorphic to $R_\mathcal{T}$.
\label{lem:RouisSerre}
\end{lem}

\proof For any compact $A$ and $B$, we have $[B,R_{\mathcal{T}}(A)] \cong [A,B]^* \cong [B,S(A)]$. By setting $B=S(A)$, we get a morphism $\nu_A: S(A) \ra R_{\mathcal{T}}(A)$. Take a cone over this morphism and denote it by $C$. Since $[B,\nu_A]$ is an isomorphism for any compact $B$, $[B,C] = 0$ for all compact $B$ and $C$ is zero. \qed

\begin{cor}
 Let $X$ be a quasi-projective scheme over a field $k$. If $D_{\perf}(X)$ possesses a weak Serre functor, it is $R_X$ and $f^!\mathcal{O}_{\Spec k}$ is perfect.
\label{cor:RouqSerre}
\end{cor}

Assume we have $k$-linear equivalences $\Phi: \mathcal{C} \ra \mathcal{C}'$ and $\Psi: \mathcal{D} \ra \mathcal{D}'$ and $k$-linear functors $F: \mathcal{C} \ra \mathcal{D}$ and $F': \mathcal{C}' \ra \mathcal{D}'$ making the diagram
\begin{center}
 \leavevmode
 \begin{xy}
  (-15,10)*+{\mathcal{C}}="b"; (15,10)*+{\mathcal{D}}="c"; {\ar@{->}^F "b";"c"}; (-15,-10)*+{\mathcal{C}'}="b'"; (15,-10)*+{\mathcal{D}'}="c'"; {\ar@{->}^{F'} "b'";"c'"}; {\ar@{->}^{\Phi} "b";"b'"}; {\ar@{->}^{\Psi} "c";"c'"};
 \end{xy}
\end{center}
commute.

\begin{lem}
 If $R_F$ exists, then so does $R_{F'}$. Moreover, $R_{F'} \circ \Phi$ is naturally isomorphic to $\Psi \circ R_F$.
\label{lem:Rouquierfunctorcommute}
\end{lem}

\proof By lemma \ref{lem:Rouqcompose}, $R_{F'} \circ \Phi \cong R_{F' \circ \Phi} \cong R_{\Psi \circ F} \cong \Psi \circ R_F$.
\qed

\begin{rmk}
 We have seen an equivalence of either the perfect derived categories or the bounded derived categories of coherent sheaves of two projective schemes induces an equivalence of the other pair of categories. By lemma \ref{lem:Rouquierfunctorcommute}, we see that these equivalences must commute with the Rouquier functors.
\label{rmk:Rouquiercommuteequiv}
\end{rmk}

\begin{prop}
 If $F$ is an exact functor between triangulated categories, then $R_F$ is also exact.
\label{prop:Roufunctriang}
\end{prop}

\proof By lemma \ref{lem:Rouquierfunctorcommute}, $R_F$ commutes with shift functors. Since both $F$ and $R_F$ commute with shifts, we have an equality $\eta_{A[1],B[1]}(\tau[1])(\sigma[1]) = \eta_{A,B}(\tau)(\sigma)$ for $\tau:B \ra R_F(A)$ and $\sigma: F(A) \ra B$.

Now we check that $R_F$ preserves triangles. Let
\begin{center}
 \leavevmode
 \begin{xy}
  (-10.5,7)*+{A}="a"; (10.5,7)*+{B}="b"; (0,-7)*+{C}="c"; {\ar@{->}^{\alpha} "a";"b"}; {\ar@{->}^{\beta} "b";"c"}; {\ar@{-->}^{\gamma} "c";"a"};
 \end{xy}
\end{center}
be an exact triangle and let $C_0$ be an object which completes the morphism $R_F(\alpha):R_F(A) \ra R_F(B)$ to an exact triangle. We now seek a morphism $\xi: C_0 \ra R_F(C)$ which makes
\begin{center}
 \leavevmode
 \begin{xy}
  (-45,10)*+{R_F(A)}="a"; (-15,10)*+{R_F(B)}="b"; (15,10)*+{C_0}="c"; (45,10)*+{R_F(A)[1]}="d"; {\ar@{->}^{R_F(\alpha)} "a";"b"}; {\ar@{->}^{\phi} "b";"c"}; {\ar@{->}^{\psi} "c";"d"};
  (-45,-10)*+{R_F(A)}="a'"; (-15,-10)*+{R_F(B)}="b'"; (15,-10)*+{R_F(C)}="c'"; (45,-10)*+{R_F(A)[1]}="d'"; {\ar@{->}_{R_F(\alpha)} "a'";"b'"}; {\ar@{->}_{R_F(\beta)} "b'";"c'"}; {\ar@{->}_{R_F(\gamma)} "c'";"d'"};
  {\ar@{->}^= "a";"a'"}; {\ar@{->}^= "b";"b'"}; {\ar@{-->}^{\xi} "c";"c'"}; {\ar@{->}^= "d";"d'"};
 \end{xy}
\end{center}
commute where the bottom row results from the application of $R_F$. If we have such a morphism, then we can apply lemma \ref{lem:diagram} to see that $\xi$ is an isomorphism. We now work on finding such a $\xi$.

Since $R_F$ is a Rouquier functor, a map $\xi \in [C_0,R_F(C)]$ is uniquely specified by an element of $\epsilon \in [F(C),C_0]^*$ where $\eta_{C,C_0}(\xi) = \epsilon$. By naturality, requiring that $\xi \circ \phi = R_F(\beta)$ is equivalent to
\begin{displaymath}
  \epsilon( \phi \circ \lambda) = \eta_{C,R_F(B)}(R_F(\beta))(\lambda) = \eta_{B,R_F(B)}(\id_{R_F(B)})(\lambda \circ F(\beta))
\end{displaymath}
for all $\lambda \in [F(C),R_F(B)]$. Let $\phi \circ [F(C),R_F(B)]$ denote the image of $[F(C),R_F(B)]$ in $[F(C),C_0]$ under post-composition with $\phi$. Consider the function
\begin{align*}
 \epsilon_1: \phi \circ [F(C),R_F(B)] & \to k \\
	   \tau & \mapsto \eta_{B,R_F(B)}(\id_{R_F(B)})(\lambda \circ F(\beta))
\end{align*}
where $\lambda$ is a choice of an element $[F(C),R_F(B)]$ satisfying $\tau = \phi \circ \lambda$. $\epsilon_1$ is well-defined. If $\lambda$ and $\lambda'$ are two elements of $[F(C),R_F(B)]$ with $\tau = \phi \circ \lambda = \phi \circ \lambda'$, then $\lambda - \lambda' = R_F(\alpha) \circ \rho$ for some $\rho \in [F(C),R_F(A)]$. By naturality,
\begin{gather*}
 \eta_{B,R_F(B)}(\id_{R_F(B)})((\lambda - \lambda') \circ F(\beta)) = \eta_{C,R_F(B)}(R_F(\beta))(\lambda - \lambda') \\ = \eta_{C,R_F(B)}(R_F(\beta))(R_F(\alpha) \circ \rho) = \eta_{C,R_F(A)}(R_F(\beta \circ \alpha))(\rho) = 0.
\end{gather*}
Any $\epsilon$ satisfying $\epsilon|_{\phi \circ [F(C),R_F(B)]} = \epsilon_1$ will make the first square commutative.

Using naturality again, requiring that $R_F(\gamma) \circ \xi = \psi$ is equivalent to 
\begin{displaymath}
 \epsilon(\mu \circ F(\gamma)) = \eta_{A[1],C_0}(\psi)(\mu) = \eta_{A[1],R_F(A)[1]}(\id_{R_F(A)[1]})(\psi \circ \mu)
\end{displaymath}
for all $\mu \in [F(A)[1],C_0]$. Let $[F(A)[1],C_0] \circ F(\gamma)$ denote the image of $[F(A)[1],C_0]$ in $[F(C),C_0]$ under pre-composition with $F(\gamma)$. Consider the function
\begin{align*}
 \epsilon_2: [F(A)[1],C_0] \circ F(\gamma) & \to k \\
	   \sigma & \mapsto \eta_{A[1],R_F(A)[1]}(\id_{R_F(A)[1]})(\psi \circ \mu)
\end{align*}
where $\mu$ is a choice of an element of $[F(A)[1],C_0]$ satisfying $\sigma = \mu \circ F(\gamma)$. $\epsilon_2$ is well-defined. If $\mu$ and $\mu'$ are two choices with $\sigma = \mu \circ F(\gamma) = \mu' \circ F(\gamma)$, then there exists a $\zeta$ in $[F(B)[1],C_0]$ with $\zeta \circ F(\alpha)[1] = \mu - \mu'$. By naturality,
\begin{gather*}
 \eta_{A[1],R_F(A)[1]}(\id_{R_F(A)[1]})(\psi \circ (\mu-\mu')) = \eta_{A[1],C_0}(\psi)(\mu-\mu') \\ = \eta_{A[1],C_0}(\psi)(\zeta \circ F(\alpha)[1]) =  \eta_{B[1],C_0}(R_F(\alpha)[1] \circ \psi)(\zeta) = 0.
\end{gather*}
Any $\epsilon$ satisfying $\epsilon|_{[F(A)[1],C_0] \circ F(\gamma)} = \epsilon_2$ will make the second square commutative. We next show that these two conditions are consistent. Namely, 
\begin{displaymath}
 \epsilon_1|_{\phi \circ [F(C),R_F(B)] \cap [F(A)[1],C_0] \circ F(\gamma)} = \epsilon_2|_{\phi \circ [F(C),R_F(B)] \cap [F(A)[1],C_0] \circ F(\gamma)}.
\end{displaymath}

Assume we have $\phi \circ \lambda = \mu \circ F(\gamma)$. We can complete the corresponding commutative square to a morphism of triangles
\begin{center}
 \leavevmode
 \begin{xy}
  (-45,10)*+{F(B)}="a"; (-15,10)*+{F(C)}="b"; (15,10)*+{F(A)[1]}="c"; (45,10)*+{F(B)[1]}="d"; {\ar@{->}^{F(\beta)} "a";"b"}; {\ar@{->}^{F(\gamma)} "b";"c"}; {\ar@{->}^{F(\alpha)[1]} "c";"d"};
  (-45,-10)*+{R_F(A)}="a'"; (-15,-10)*+{R_F(B)}="b'"; (15,-10)*+{C_0}="c'"; (45,-10)*+{R_F(A)[1]}="d'"; {\ar@{->}_{R_F(\alpha)} "a'";"b'"}; {\ar@{->}_{\phi} "b'";"c'"}; {\ar@{->}_{\psi} "c'";"d'"};
  {\ar@{-->}^{\nu} "a";"a'"}; {\ar@{->}^{\lambda} "b";"b'"}; {\ar@{->}^{\mu} "c";"c'"}; {\ar@{-->}^{\nu[1]} "d";"d'"};
 \end{xy}
\end{center}
which gives the relations $\lambda \circ F(\beta) = R_F(\alpha) \circ \nu$ and $\psi \circ \mu = \nu[1] \circ F(\alpha)[1]$. Using the naturality of $\eta$ and its commutation with shifts, we have the following equalities
\begin{gather*}
 \eta_{B,R_F(B)}(\id_{R_F(B)})( \lambda \circ F(\beta)) = \eta_{B,R_F(B)}(\id_{R_F(B)})( R_F(\alpha) \circ \nu ) = \\ \eta_{B,R_F(A)}(R_F(\alpha))(\nu) = \eta_{A,R_F(A)}(\id_{R_F(A)})( \nu \circ F(\alpha) ) = \\ \eta_{A[1],R_F(A)[1]}(\id_{R_F(A)[1]})( \nu[1] \circ F(\alpha)[1]) = \eta_{A[1],R_F(A)[1]}(\id_{R_F(A)[1]})(\psi \circ \mu).
\end{gather*}

Thus, we have a well-defined linear map, which we shall denote by
\begin{displaymath}
 \epsilon: \phi \circ [F(C),R_F(B)] \cup [F(A)[1],C_0] \circ F(\gamma) \to k,
\end{displaymath}
whose restriction to $\phi \circ [F(C),R_F(B)]$ is $\epsilon_1$ and whose restriction to $[F(A)[1],C_0] \circ F(\gamma)$ is $\epsilon_2$. As we are working with vector spaces over $k$, we can choosing a splitting over $k$,
\begin{displaymath}
 [F(C),C_0] \cong \left(\phi \circ [F(C),R_F(B)] \cup [F(A)[1],C_0] \circ F(\gamma) \right) \oplus W.
\end{displaymath}
To define our desired $\epsilon$, it is sufficient to specify it on each summand. We define it to be $\epsilon$, as above, on $\phi \circ [F(C),R_F(B)] \cup [F(A)[1],C_0] \circ F(\gamma)$ and we take it to be zero (for specificity) on $W$. This gives the desired $\epsilon$ and, consequently, the desired $\xi$.  \qed

\section{An application: reconstruction for projective Gorenstein varieties}
\label{sec:recon}
 
For this section, a variety is a reduced and irreducible scheme of finite type over a field $k$. In \cite{BO01}, Bondal and Orlov prove the following reconstruction result.
\begin{thm}
 Let $X$ be a smooth projective variety over a field $k$ with ample or anti-ample canonical bundle. Assume there is another smooth variety, $Y$, and an exact equivalence $D^b_{\coh}(X) \cong D^b_{\coh}(Y)$. Then, $X$ is isomorphic to $Y$.
\end{thm}

In this section, we use the results of the previous sections to augment the original argument of Bondal and Orlov and prove the following extension. 

\begin{prop}
 Let $X$ be a projective Gorenstein variety over a field $k$ with ample or anti-ample canonical bundle. Assume that $Y$ is a projective variety over $k$ and there is an exact equivalence between $D_{\perf}(X)$ and $D_{\perf}(Y)$. Then, $X$ is isomorphic to $Y$.
\label{prop:reconstructionvariety}
\end{prop}

An immediate corollary to proposition \ref{prop:reconstructionvariety}, thanks to lemma \ref{lem:misclfdualresults}, is the following.

\begin{cor}
 Let $X$ be a projective Gorenstein variety over a field $k$ with ample or anti-ample canonical bundle. Assume that $Y$ is a projective variety over $k$ and there is an exact equivalence between $D^b_{\coh}(X)$ and $D^b_{\coh}(Y)$. Then, $X$ is isomorphic to $Y$.
\end{cor}

Recall that a Noetherian scheme, $X$, is called \textbf{Gorenstein} if, for each point $x$ of $X$, the local ring $\mathcal O_{x,X}$ has finite injective dimension as a module over itself. We have the following more useful characterization from proposition 9.3 of \cite{RD}.

\begin{prop}
 Let $X$ be a connected projective scheme. $X$ is Gorenstein if and only $f^!\mathcal O_{\operatorname{Spec} k}$ is quasi-isomorphic to a line bundle concentrated in a single degree, $-\dim X$, where $f: X \to \operatorname{Spec} k$ is the structure morphism.
\end{prop} 

\begin{defn}
 A scheme $X$ over a field $k$ is \textbf{categorically Gorenstein} if $D_{\perf}(X)$ has a Serre functor.
\end{defn}

\begin{lem}
 If $X$ is a projective variety, $X$ is categorically Gorenstein if and only if it is Gorenstein. In particular, $f^!\mathcal{O}_{\Spec k}$ is a line bundle concentrated in degree $-\dim X$.
\label{lem:catGor=Gor}
\end{lem}

\proof If $X$ is Gorenstein, then $f^!\mathcal{O}_{\Spec k}$ is a line bundle concentrated in $-\dim X$ which implies that $R_X$ is a Serre functor.

If $X$ is categorically Gorenstein, we know, by corollary \ref{cor:RouqSerre}, $R_X$ must be the Serre functor. In particular, $f^!\mathcal{O}_{\Spec k}$ must be perfect. $R_X$ has a right adjoint given by tensoring with the dual complex. As $R_X$ is an autoequivalence, the adjoint must be the inverse.

Assume we are given two perfect complexes, $\mathcal E$ and $\mathcal E^{-1}$, so that the functors $-\otimes \mathcal E$ and $- \otimes \mathcal E^{-1}$ are inverses. Note that $\mathcal E^{-1}$ must be $\mathcal E^{\vee}=\mathcal{H}om(E,\mathcal O_X)$ as $-\otimes \mathcal E^{\vee}$ is left and right adjoint to $- \otimes \mathcal E$ on $D(X)$. By proposition \ref{prop:exisuniq} and lemma \ref{lem:misclfdualresults}, if $-\otimes \mathcal E$ is an autoequivalence of $D_{\perf}(X)$, then $-\otimes \mathcal E^{-1}$ is an autoequivalence of $D^b_{\coh}(X)$. This, in turn, implies that $-\otimes \mathcal E$ is an autoequivalence on $D^b_{\coh}(X)$.

Take a point, $x$, of $X$ and pull back to the local ring $\mathcal{O}_{x,X}$. Let $\mathcal E_x$ and $\mathcal E_x^{-1}$ denote the pullbacks of $\mathcal E$ and $\mathcal E^{-1}$ respectively. Up to quasi-isomorphism, we can assume that $\mathcal E_x$ and $\mathcal E_x^{-1}$ are both minimal complexes. We must have
\begin{displaymath}
 \mathcal E_x \otimes_{\mathcal{O}_{x,X}} \mathcal E_x^{-1} \otimes_{\mathcal O_{x,X}} k(x) \cong k(x).
\end{displaymath}
Rewriting, we see
\begin{displaymath}
 \left(\mathcal E_x \otimes_{\mathcal{O}_{x,X}} k(x) \right) \otimes_{k(x)} \left(\mathcal E_x^{-1} \otimes k(x)\right) \cong k(x).
\end{displaymath}
The complex, $\mathcal E_x \otimes_{\mathcal{O}_{x,X}} k(x)$, has zero differential and its rank in degree $i$ is the same as the rank of $\mathcal E_x$ in degree $i$. The same holds from $\mathcal E_x^{-1} \otimes k(x)$. The only way to tensor these two graded vector spaces together and get $k(x)$ in degree zero is for $\mathcal E_x \otimes_{\mathcal{O}_{x,X}} k(x)$ to be isomorphic to $k(x)[i]$ and $\mathcal E_x^{-1} \otimes_{\mathcal{O}_{x,X}} k(x)$ to be isomorphic to $k(x)[-i]$ for some $i$. Thus, $\mathcal E_x$ is quasi-isomorphic to $\mathcal{O}_{x,X}[i]$. This holds on some neighborhood of $x$ and $i$ is locally-constant as a function of $x$. Since $X$ is connected, $i$ is constant and $\mathcal E_x$ is the shift of a line bundle.

Consequently, $f^! \mathcal O_{\Spec k}$ is the shift of a line bundle. The shift must be the dimension of $X$, as can be checked by restricting to a closed point. \qed

%

If $X$ is Gorenstein, we have a Serre functor $S$ on $D_{\perf}(X)$. $S^{\vee}$, the right pseudo-adjoint to $S$, induces an autoequivalence $D^b_{\coh}(X) \ra D^b_{\coh}(X)$, which must be $- \otimes \omega_X^{-1}[-\dim X]$ by proposition \ref{prop:exisuniq}.  We shall use this to characterize points a'la \cite{BO01}.

\begin{defn}
 A \textbf{point functor} $P$ of codimension $d$ of $D_{\perf}(X)$ is a locally-finite cohomological functor on $D_{\perf}(X)$, which satisfies the following conditions
\begin{enumerate}
 \item $S^{\vee}(P) \cong P[-d]$.
 \item $[P,P[l]]$ is zero for $l < 0$.
 \item $[P,P] \cong k(P)$ a finite field extension of $k$.
\end{enumerate}
\end{defn}

\begin{lem}
 Let $X$ be a projective Gorenstein variety. Any point object must have codimension $\dim X$.
\end{lem}

\proof By uniqueness of pseudo-adjoints, proposition \ref{prop:exisuniq},
\begin{displaymath}
 S^{\vee} \cong \omega_X^{-1}[-\dim X] \otimes -.
\end{displaymath}
By theorem \ref{thm:lfcohom}, $P$ is represented by a bounded complex of coherent sheaves which we also denote by $P$. Let $P$ have codimension $d$. Since $S^{\vee} (P) \cong P[-d]$, we know that $\omega_X^{-1} \otimes P$ is quasi-isomorphic to $P[-d+\dim X]$. Let $\mathcal{H}^i$ denote the cohomology sheaves of $P$. Since $P$ has bounded cohomology and $\omega_X^{-1} \otimes \mathcal{H}^i \cong \mathcal{H}^{i-d+\dim X}$, either $P$ is quasi-isomorphic to zero (and is not a point object) or $d=\dim X$. \qed

\begin{lem}
 Let $X$ be a projective Gorenstein variety over a field $k$ with ample or anti-ample canonical bundle. Then, an object of $P$ of $D^b_{\coh}(X)$ is a point object if and only if $P$ is isomorphic to $\mathcal{O}_p[r]$ for some closed point $p$ of $X$.
\label{lem:point}
\end{lem}

\proof We shall follow the proof in \cite{Huy06}. Note that shifts of points are point objects of codimension $\dim X$.
Since $S^{\vee} (P) \cong P[-\dim X]$, we know that $\omega_X^{-1} \otimes P$ is quasi-isomorphic to $P$. Let $\mathcal{H}^i$ denote the cohomology sheaves of $P$. Because $\omega_X^{-1}$ is ample (or anti-ample) and $\omega_X^{-1} \otimes \mathcal{H}^i \cong \mathcal{H}^i$, $\mathcal{H}^i$ has zero dimensional support.

We can resolve $P$ using direct sums of injective sheaves each of whose support is contained within an irreducible component of the support of the cohomology sheaves of $P$. Since any map between two quasi-coherent sheaves with disjoint support is zero, we get a splitting of $P$ into complexes supported at single points. Since $[P,P]$ is a field, all summands but one must be quasi-isomorphic to zero.

Assume $P$ has cohomology supported only at a single point of $X$. Let $m_0$ be the minimal $i$ so that $\mathcal{H}^i$ is nonzero and $m_1$ the maximal $i$ so that $\mathcal{H}^i$ is nonzero. By truncating, we can assume that $P$ is zero outside $[m_0,m_1]$.  And there are morphisms $\mathcal{H}^{m_0} \ra P[m_0]$ and $P[m_1] \ra \mathcal{H}^{m_1}$.  For each $\mathcal{H}^{m_0}$ and $\mathcal{H}^{m_1}$, there are nonzero maps in and out of $k(x)$. Composing these gives a nontrivial element of $[P[m_1],P[m_0]]$. Thus, $m_0$ equals $m_1$ and $P$ is a shift of a coherent sheaf. If the length of $P$ is greater than one we can project down the composition series to get a non-invertible endomorphism. Thus, $P$ is simply $\mathcal{O}_p[r]$ for some closed point $p$ and some $r \in \Z$. \qed

\begin{defn}
 An object $L$ of $D_{\perf}(X)$ is an \textbf{locally-free object} if, for any point object, $P$, there exists a $t \in \Z$  and $n > 0$ so that
\begin{enumerate}
 \item $[L,P[t]] \cong k(P)^n$ and
 \item $[L,P[i]] = 0$ for $i \not = t$.
\end{enumerate}
\label{def:locfree}
\end{defn}

\begin{lem}
 Let $X$ be a projective variety. A perfect object, $L$, satisfies the two conditions of definition \ref{def:locfree} for all shifts of closed points if and only if $L$ is isomorphic to a shift of a locally-free sheaf. In particular, if all the point objects of $X$ are isomorphic to shifts of points, then $L$ is an locally-free object if and only if $L$ is a shift of a locally-free coherent sheaf.
\label{lem:invert}
\end{lem}



\proof Take a closed point $x \in X$. Using the minimal resolution of $L_x$ in an argument similar to the proof of lemma \ref{lem:catGor=Gor}, we conclude that $L_x$ must be quasi-isomorphic to a free module of rank $n$ concentrated in a fixed degree. As $X$ is connected, $L$ must have cohomology only in a single degree, with that cohomology sheaf being locally-free. \qed


Now we move onto the proof of proposition \ref{prop:reconstructionvariety}. It will be accomplished through a sequence of lemmas. Assume that $X$ and $Y$ satisfy the hypotheses of proposition \ref{prop:reconstructionvariety}. We also assume, for definiteness, that $\omega_X$ is ample. The case where $\omega_X^{-1}$ is ample follows from a similar argument. $D$ will be a placeholder for $D_{\perf}(X)$ or $D_{\perf}(Y)$ depending on the context. Since they are equivalent, this should cause no confusion.

\begin{lem}
 All point functors of $Y$ are shifts of points.
\end{lem}

\proof Let us denote the category of point objects of $Z$ as $\mathcal{P}(Z)$. From lemma \ref{lem:point}, we know there is a bijection between shifts of closed points of $X$ and objects of $\mathcal{P}(X)$. From our assumption, we have a equivalence between $\mathcal{P}(X)$ and $\mathcal{P}(Y)$. By lemma \ref{lem:catGor=Gor}, $Y$ is Gorenstein and, thus, there in inclusion of the shifts of points on $Y$ into $\mathcal{P}(Y)$. The category $\mathcal{P}(X)$ satisfies the following condition: if $P$ and $Q$ are objects, then either $P$ is isomorphic to a shift of $Q$ or $[P,Q[j]]$ is zero for all integers $j$. If $N$ is nonzero bounded complex of coherent sheaves, then let $m$ be the largest integer for which the $m$-th cohomology sheaf is nonzero. There must be a nonzero map in $[N,\mathcal{O}_y[-m]]$ for some $y$. Thus, $N$ must be a shift of a point. So all point objects of $Y$ are shifts of points. \qed

We can then apply lemma \ref{lem:invert}.

\begin{lem}
 All locally-free objects of $D_{\perf}(Y)$ are shifts of locally-free coherent sheaves.
\end{lem}

\begin{lem}
 The sets of closed points of $X$ and $Y$, with their Zariski topologies, are homeomorphic.
\end{lem}

\proof Choose an invertible object $L_0$ corresponding to an invertible sheaf on $X$. By shifting, we can assume our equivalence takes $L_0$ to a complex quasi-isomorphic to an invertible sheaf on $Y$. Let us denote the image by $L_0$ also. Now the set of point objects $P$ so that $[L_0,P]$ is $k(P)$ is in bijection with the set of closed points of $X$. Denote this set by $p_D$. Similarly, $p_D$ is in bijection with the closed points of $Y$. This gives us a bijection between the closed points of $X$ and $Y$. 

Let $l_D$ denote the set of locally-free objects $L$ in $D$ so that $[L,P]$ is isomorphic to $k(P)^n$ for some $n$ and for all $P$ in $p_D$. $l_D$ is in bijection with the set of locally-free sheaves on $X$ and the set of locally-free sheaves on $Y$. For $\alpha$ in $[L,L']$, with $L,L'$ in $l_D$, let $U_{\alpha}$ denote the set of $P$ in $p_D$ so that the induced map
\begin{displaymath}
 - \circ \alpha : [L',P] \ra [L,P]
\end{displaymath}
is nonzero. Since $X$ and $Y$ possess enough locally-free coherent sheaves, we know the open sets $U_{\alpha}$ in $X$ and $Y$ form a basis for the topologies of $X$ and $Y$, see \cite{Ill71}. Thus, our identification of points is a homeomorphism. \qed

Therefore, the dimensions of $X$ and $Y$ must coincide. Record this common dimension as $d$.

\begin{lem}
 $\omega_Y$ is ample.
\end{lem}

\proof Let $L$ be a line bundle on an algebraic variety $V$. $U_{\alpha}$ for $\alpha$ in $[L^{\otimes i},L^{\otimes j}]$ form a basis for the topology of $V$ if and only if $L$ is ample \cite{Ill71}. We see that our equivalence takes ample invertible sheaves to shifts of ample invertible sheaves.

We can twist our equivalence by an invertible sheaf and shift it to assume that the structure sheaf of $X$ is sent to the structure sheaf of $Y$. Then, from the naturality of the Serre functors, $\omega_X[d]$ is sent to $\omega_Y[d]$. \qed

\begin{lem}
 The graded rings $\bigoplus_{n\in\Z} H^0(\mathcal{O}_X,\omega^n_X)$ and $\bigoplus_{n\in\Z} H^0(\mathcal{O}_Y,\omega^n_Y)$ are isomorphic. Consequently, $X$ and $Y$ are isomorphic.
\end{lem}

\proof Set $L_i$ equal to $S^iL_0[-di]$ where $L_0$ is as chosen before. For each pair $(i,j)$, we have natural isomorphisms
\begin{gather*}
 [L_i.L_j] \cong [S^iL_0[-di],S^jL_0[-dj]] \cong [L_0,S^{j-i}L_0[-d(j-i)]] \cong [L_0,L_{j-i}]
\end{gather*}
This provides the structure of a graded ring for $A = \bigoplus_{l=-\infty}^{\infty} \Hom_D(L_0,L_l)$. But, $A$ is isomorphic to $\bigoplus_{l\in\Z} H^0(X,\omega_X^l)$ and $\bigoplus_{l\in\Z} H^0(Y,\omega_Y^l)$. Since both $\omega_X$ and $\omega_Y$ are either ample, we can take Proj of the non-negative portion to give $X \cong Y$. \qed

\begin{rmk}
 Reconstruction can be accomplished with less hypotheses on $X$ and $Y$. We can replace ampleness of $f^!\mathcal O_{\Spec k}[-\dim X]$ by the condition that $f^!\mathcal O_{\Spec k}$ is perfect and the smallest thick subcategory of $D_{\perf}(X)$ containing all positive tensor powers of $f^!\mathcal O_{\Spec k}$ is $D_{\perf}(X)$ itself. The argument above is no longer sufficient under this assumption. However, one can used a modified argument of Rouquier, see \cite{Rou02}. 
\end{rmk}

One can also carry forth the arguments used in \cite{BO01} to prove the following.

\begin{prop}
 Let $X$ be a projective Gorenstein variety with ample or anti-ample canonical bundle. Then, the group of auto-equivalences of $D_{\perf}(X)$ (and $D^b_{\coh}(X)$) is isomorphic to $\Aut(X) \rtimes (\Pic(X) \times \Z)$.
\end{prop}

\proof (Sketch) We only provide a sketch omitting some details from the arguments in \cite{BO01}. Assume for definiteness that $\omega_X$ is ample. We know by lemma \ref{lem:misclfdualresults} that the group of autoequivalences of $D_{\perf}(X)$ and $D^b_{\coh}(X)$ coincide. So, we can assume we have an autoequivalence $\Phi$ of $D^b_{\coh}(X)$ whose restriction to $D_{\perf}(X)$ is also an autoequivalence. $\Phi$ must commute with the Serre functor on $D_{\perf}(X)$ so $\Phi$ must also commute with $S^{\vee}$ on $D^b_{\coh}(X)$ by proposition \ref{prop:exisuniq}. Now we proceed exactly as in \cite{BO01}.

As point functors and locally-free objects are defined in a purely categorical manner, $\Phi$ must take point functors to point functors and locally-free objects to locally-free objects. Using the argument above, we get that, after possibly shifting and tensoring with a line bundle, $\Phi$ induces an action on $\bigoplus_{n\in\Z} H^0(\mathcal{O}_X,\omega^n_X)$. This gives an element $\phi \in \op{Aut}(X)$. Replacing $\Phi$ by $\Phi \circ \left(\phi^{-1}\right)^*$ we can assume that $\Phi$ acts trivially on $\op{Proj}\left(\bigoplus_{n\geq 0} H^0(\mathcal{O}_X,\omega^n_X)\right)$. After appropriately rescaling $\Phi$, this gives a natural transformation $\op{Id} \to \Phi$ defined on all coherent sheaves and an isomorphism on the subcategory consisting of powers of the canonical bundle. Applying proposition A.3 of \cite{BO01} shows that $\Phi$ must be isomorphic to the identity functor on $D^b_{\coh}(X)$. \qed

\end{document}